\documentclass[12pt]{amsart}

\usepackage[mathscr]{eucal}
\usepackage{amsthm,amscd,amssymb}
\usepackage[all]{xy}
\xyoption{arc} \xyoption{rotate}

\newcommand{\non}{\nonumber}
\newcommand{\wt}{\widetilde}

\newcommand{\la}{\lambda}

\newcommand{\ga}{\gamma}
\newcommand{\Ga}{\Gamma}

\newcommand{\ts}{\,}
\newcommand{\tss}{\hspace{1pt}}

\newcommand{\Y}{ {\rm Y}}

\newcommand{\End}{ {\rm End}}
\newcommand{\C}{\mathbb{C}\tss}
\newcommand{\Z}{\mathbb{Z}\tss}

\newcommand{\gl}{\mathfrak{gl}}
\newcommand{\g}{\mathfrak{g}}
\newcommand{\sll}{\mathfrak{sl}}

\newcommand{\cdet}{ {\rm cdet}\ts}
\newcommand{\gr}{ {\rm gr}}
\newcommand{\ad}{ {\rm ad}}
\newcommand{\sgn}{ {\rm sgn}}

\newcommand{\dmo}{\text{-mod}}

\newtheorem{thm}{Theorem}[section]
\newtheorem{prop}[thm]{Proposition}
\newtheorem{cor}[thm]{Corollary}
\newtheorem{lem}[thm]{Lemma}

\theoremstyle{definition}
\newtheorem{defin}[thm]{Definition}

\theoremstyle{remark}
\newtheorem{remark}[thm]{Remark}
\newtheorem{example}[thm]{Example}

\newcommand{\bth}{\begin{thm}}
\renewcommand{\eth}{\end{thm}}
\newcommand{\bpr}{\begin{prop}}
\newcommand{\epr}{\end{prop}}
\newcommand{\ble}{\begin{lem}}
\newcommand{\ele}{\end{lem}}
\newcommand{\bco}{\begin{cor}}
\newcommand{\eco}{\end{cor}}
\newcommand{\bde}{\begin{defin}}
\newcommand{\ede}{\end{defin}}
\newcommand{\bex}{\begin{example}}
\newcommand{\eex}{\end{example}}
\newcommand{\bre}{\begin{remark}}
\newcommand{\ere}{\end{remark}}
\newcommand{\bal}{\begin{aligned}}
\newcommand{\eal}{\end{aligned}}
\newcommand{\beq}{\begin{equation}}
\newcommand{\ben}{\begin{equation*}}

\def\beql#1{\begin{equation}\label{#1}}

\newcommand\dcl\DeclareMathOperator

\dcl\Sp{Specm} \dcl\St{St} \dcl\cfs{cfs} \dcl\supp{supp}
\dcl\Ker{Ker} \dcl\Hom{Hom} \dcl\Ext{Ext} \dcl\Ann{Ann} \dcl\Ob{Ob}
\dcl\im{Im} \dcl\mo{mod} \dcl\rank{rank} \dcl\M{M} \dcl\Specm{Specm}
\dcl\Aut{Aut} \dcl\lm{lm} \dcl\lc{lc} \dcl\lt{lt} \dcl\Gal{Gal}

\newcommand{\mL}{\mathcal L}

\newcommand{\bm}{\mathbf m}
\newcommand{\bn}{\mathbf n}

\newcommand\osu[2]{\overset{#2}{\underset{#1}\oplus}\,}

\newcounter{numberofremark}
\DeclareMathOperator\gro{growth}

\dcl\gkdim{GKdim}

\newcommand\vi{\varphi}

\newcommand\de{{\delta}}

\newcommand\bZ{{\mathbb Z}}

\newcommand\A{{{\mathscr A}}}

\newcommand\cA{{\mathscr A}}

\newcommand\cZ{{\mathscr Z}}

\renewcommand{\k}{\Bbbk}

\newcommand\ds{\displaystyle}

\newcommand\myto{\longrightarrow}
\newcommand\pa{\partial}

\newcommand\cM{\mathscr M}

\newcommand{\cond}[2]{
\ifthenelse {\equal{\commenttag} {detailed}} {#2$\maltese$}
 {#1}
 }

\dcl\Der{Der} \dcl\Inn{Inn} \dcl\mcd{mcd} \dcl\GK{GK} \dcl\Ass{Ass}
\dcl\Spec{Spec}

\newcommand\cO{\mathscr{O}}

\newcommand\bF{\mathbb F}
\newcommand\bN{\mathbb N}

\newcommand{\eeq}{\end{equation}}
\newcommand{\een}{\end{equation*}}

\def\bm{\mathbf m}

\usepackage{color}
\renewcommand\k{\Bbbk}

\begin{document}

\title[Gelfand-Kirillov conjecture]
{GELFAND-KIRILLOV CONJECTURE  \\[0.4em]
AND GELFAND-TSETLIN MODULES FOR\\[0.4em]
FINITE $W$-ALGEBRAS}

\author{Vyacheslav Futorny}
\address{Institute of Mathematics and Statistics\\
University of S\~ao Paulo\\
Caixa Postal 66281- CEP 05315-970\\
S\~ao Paulo, Brazil} \email{futorny@ime.usp.br}

\author{Alexander Molev}
\address{
School of Mathematics and Statistics\\
University of Sydney, NSW 2006, Australia}
\email{alexm@maths.usyd.edu.au}

\author{Serge Ovsienko}
\address{
Faculty of Mechanics and Mathematics\\
Kiev Taras Shevchenko University\\
Vla\-di\-mir\-skaya 64, 00133, Kiev, Ukraine}
\email{ovsyenko@univ.kiev.ua}

\begin{abstract}
We address two problems regarding the structure and representation
theory of finite $W$-algebras associated with the general linear
Lie algebras. Finite $W$-algebras can be defined either via the
Whittaker modules of Kostant or, equivalently, by the  quantum
Hamiltonian reduction. Our first main result is a proof of the
Gelfand-Kirillov conjecture for the skew fields of fractions of
the finite $W$-algebras. The second main result is a
parametrization of finite families of irreducible Gelfand-Tsetlin
modules by the characters of the Gelfand-Tsetlin subalgebra. As a
corollary, we obtain a complete classification of generic
irreducible Gelfand-Tsetlin modules for the finite $W$-algebras.
\end{abstract}

\vspace{0,2cm}

\maketitle


\newpage

\section {Introduction}
\setcounter{equation}{0}

The concept of finite $W$-algebras goes back to the original paper
of Kostant \cite{Ko} dealing with the study of Whittaker modules
and to its generalization by Lynch \cite{L}. An alternative
construction of $W$-algebras can be given via the  quantum
Hamiltonian reduction which goes back to the works of Feigin and
Frenkel \cite{FF}, Kac, Roan and Wakimoto \cite{KRW}, Kac and
Wakimoto \cite{KW} and  De Sole and Kac \cite{SK}. It was shown by
D'Andrea, De Concini, De Sole, Heluani and Kac \cite[Appendix]{SK}
and by Arakawa~\cite{A} that both definitions of finite
$W$-algebras are equivalent.

Let $\g=\gl_m$ denote the general linear Lie algebra over
an algebraically closed field
${\Bbbk}$ of characteristic $0$ which will be fixed
throughout the paper.
A finite $W$-algebra can be associated to a fixed nilpotent element
$f\in \g$ as follows. A $\Z$-grading $\g=\oplus^{}_{j\in \Z}\g_j$
is called a \emph{good grading} for $f$ if
$f\in \g_2$ and the linear map
\ben
\ad\ts f:\g_j\to\g_{j+2}
\een
is injective
for $j\leqslant -1$ and surjective for
$j\geqslant -1$. A complete classification of good gradings
for simple Lie algebras was given by
Elashvili and Kac \cite{EK}.
A non-degenerate invariant symmetric bilinear form $(.\ts,.)$ on $\g$
induces a non-degenerate skew-symmetric form on
$\g_{-1}$
defined by $\langle x,y\rangle=([x,y],\ts f)$.
Let $\mathcal{I}\subset \g_{-1}$ be a
maximal isotropic subspace and
set $\mathfrak{t}=\bigoplus_{j\leqslant
-2}\g_j \oplus\mathcal{I}$. Now let
$\chi:U(\mathfrak{t})\rightarrow \C$
be the one-dimensional representation such that
$x \mapsto (x,f)$
for any $x\in \mathfrak{t}$. Set $I_{\chi}=\Ker \chi$ and
$Q_{\chi}=U(\g)/U(\g)I_{\chi}$. The corresponding finite
$W$-algebra is
defined by
\ben
W(\chi)=\End_{U(\g)}(Q_{\chi})^{op}.
\een
If the grading on $\g$ is \emph{even}, i.e. $\g_j=0$ for all odd $j$,
then $W(\chi)$ is isomorphic to the subalgebra of $\mathfrak{t}$-twisted
invariants in $U(\mathfrak{p})$ for the parabolic subalgebra
$\mathfrak{p}=\oplus^{}_{j\geqslant 0}\g_j$. Note that by the results of
Elashvili and Kac \cite{EK}, it is sufficient to consider only even
good gradings.

The growing interest to the theory of finite $W$-algebras is due, on
the one hand, to their geometric realizations
as quantizations of the Slodowy slices (see Premet~\cite{P}
and Gan and Ginzburg~\cite{GG}),
and, on the other hand, to their close
connections with the Yangian theory which was
originally observed by Ragoucy and Sorba~\cite{RS}
and developed in full generality by
Brundan and Kleshchev~\cite{BK1}.
The latter results may well be regarded
as a substantial step forward in
understanding the structure of the finite $W$-algebras
associated to $\gl_m$.
These algebras turn out to be isomorphic to
certain quotients of the \emph{shifted Yangians}, which
provides their presentations in terms of generators
and defining relations and thus opens the way for
developing the representation theory for the finite $W$-algebras;
see \cite{BK2}.

In more detail, following \cite{EK},
consider a {\it pyramid\/}
$\pi$ which is a unimodal sequence
$(q_1,q_2,\dots,q_l)$ of positive integers with
$q_1\leqslant \dots\leqslant q_k$ and
$q_{k+1}\geqslant\dots \geqslant q_l$
for some $0\leqslant k\leqslant l$.
Such a pyramid can be visualized as the diagram
of bricks (unit squares) which consists of $q_1$ bricks
stacked in the first (leftmost) column, $q_2$ bricks
stacked in the second column, etc.
The pyramid $\pi$ defines the tuple $(p_1,\dots,p_n)$ of its row
lengths, where $p_i$ is the number of bricks in the $i$th row of the
pyramid, so that $1\leqslant p_1\leqslant \dots\leqslant p_n$.
The figure illustrates the pyramid with the columns (1,3,4,2,1) and
rows (1,2,3,5):
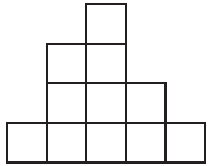
\ben {\begin{picture}(90, 60) \put(0,0){\line(1,0){75}}
\put(0,15){\line(1,0){75}} \put(15,30){\line(1,0){45}}
\put(15,45){\line(1,0){30}} \put(30,60){\line(1,0){15}}
\put(0,0){\line(0,1){15}} \put(15,0){\line(0,1){45}}
\put(30,0){\line(0,1){60}} \put(45,0){\line(0,1){60}}
\put(60,0){\line(0,1){30}} \put(75,0){\line(0,1){15}}
\end{picture}}
\een

If
the total number of bricks in the pyramid $\pi$ is $m$, then the
finite $W$-algebra $W(\pi)$ associated to $\gl_m$ corresponds to the
nilpotent matrix $f\in \gl_m$ of Jordan type $(p_1,\dots,p_n)$; see
Section~\ref{sec:sy} for the precise definition and the relationship
of $W(\pi)$ with the shifted Yangian. One of surprising consequences
of the results of \cite{BK1} is that the isomorphism class of
$W(\pi)$ depends only on the sequence of row lengths
$(p_1,\dots,p_n)$ of $\pi$. Therefore, we may assume
without restricting generality that the rows of $\pi$
are left-justified.

The first problem we address in this paper is the
\emph{Gelfand-Kirillov conjecture} for the algebras $W(\pi)$. The
original conjecture states that the universal enveloping algebra
of an algebraic Lie algebra over an algebraically closed field is
"birationally" equivalent to some Weyl algebra over a purely
transcendental extension of $\k$, i.e. its skew field of fractions
is a Weyl field. The conjecture was settled in the original paper by
Gelfand and Kirillov~\cite{GK1} for nilpotent Lie algebras, and for
$\gl_m$ and $\sll_m$; see also \cite{GK2}, where its weaker form was
proved. For solvable Lie algebras the conjecture was settled by
Borho, Gabriel and Rentschler \cite{BGR}, Joseph \cite{Jo} and
McConnell \cite{Mc}. Some mixed cases were considered by
Nghiem~\cite{Ng}, while Alev, Ooms and Van den Bergh~\cite{AOV1}
proved the conjecture for all Lie algebras of dimension at most
eight. On the other hand, counterexamples to the conjecture are
known for certain semi-direct products; see e.g. \cite{AOV2}. We
refer the reader to the book by Brown and Goodearl~\cite{BG} and
references therein for generalizations of the Gelfand-Kirillov
conjecture for quantized enveloping algebras.

For an associative algebra $A$ we  denote by $D(A)$ its skew field
of fractions, if it exists. Let $A_k$ be the $k$-th Weyl algebra
over $\k$ and $D_k=D(A_k)$ its skew field of fractions.
 Let $\mathcal{F}$ be a pure transcendental
extension of $\k$ of degree $m$ and let $A_k(\mathcal{F})$ be the
$k$-th Weyl algebra over $\mathcal F$. Denote by $D_{k,m}$ the skew
field of fractions of $A_k(\mathcal{F})$.

\medskip

\emph{Gelfand-Kirillov problem for $W(\pi)$: Does
$D(W(\pi))\simeq D_{k,m}$ for some $k,\ts m$?}

\medskip

Our first main result is a positive solution of
this problem.

\medskip

\noindent{\bf Theorem I.} The Gelfand-Kirillov conjecture holds for
$W(\pi)$:
$$D(W(\pi))\simeq D_{k,m},$$
where
$k=\sum_{i=1}^l q_i(q_i-1)/2$ and $m=q_1+\ldots+q_l$.

\medskip

Note that $m$ is the number of bricks in the pyramid $\pi$, while
$k$ can be interpreted as the sum of all leg lengths of the bricks.
Hence, $k$ and $m$ can be expressed in terms of the rows as
$k=(n-1)\ts p_1+\ldots +p_{n-1}$ and $m=p_1+\ldots+p_n$.
In the case of the one-column pyramid $(1,\dots,1)$ of height $m$
we recover the original result of \cite{GK1} for $\gl_m$.
One of the key points in the proof of Theorem I is a positive
solution of the \emph{noncommutative Noether problem} for the symmetric
group $S_k$:

\medskip
\emph{Noncommutative Noether problem for $S_k$: Does
$D_k^{S_k}\simeq D_k$?}

\medskip

Here $S_k$ acts naturally on $A_k$ and on $D_k$ by simultaneous
permutations of variables and derivations.

The second problem that we address in this paper is the
classification problem of irreducible Gelfand-Tsetlin modules
(sometimes also called Harish-Chandra modules) for finite
$W$-algebras with respect to the Gelfand-Tsetlin subalgebra. Given
a pyramid $\pi$ with the left-justified
rows $(p_1,\dots,p_n)$, for each $k\in\{1,\dots,n\}$ we let $\pi_k$
denote the pyramid $(p_1,\dots,p_k)$. We have the
chain of natural subalgebras
\beql{chainw}
W(\pi_1)\subset
W(\pi_2) \subset\dots\subset W(\pi_n)=W(\pi).
\eeq
Denote by
$\Gamma$ the (commutative) subalgebra of $W(\pi)$ generated by the
centers of the subalgebras $W(\pi_k)$ for $k=1,\dots,n$. Note that
the structure of the center of the algebra $W(\pi)$ is described
in \cite[Theorem~6.10]{BK2}. Following the terminology of that
paper, we call $\Gamma$ the \textit{Gelfand--Tsetlin subalgebra}
of $W(\pi)$.

A finitely generated module $M$ over $W(\pi)$ is called a {\em
Gelfand-Tsetlin module\/} (with respect to $\Gamma$) if \ben
M=\underset{{\bm} \in \Specm {\Gamma}}{\bigoplus}M({\bm}) \een as
a $\Ga$-module, where \ben M({\bm}) \ = \ \{ x\in M\ | \ {\bm}^k x
=0\quad \text{for some}\quad k\geqslant 0\} \een and $\Specm \Ga$
denotes the set of maximal ideals of $\Ga$. In the case of the
one-column pyramids $\pi$ this reduces to the definition of the
Gelfand--Tsetlin modules for $\gl_m$ \cite{dfo:gz}. Note also that
the {\it admissible} $W(\pi)$-modules of \cite{BK2} are
Gelfand-Tsetlin modules.

An irreducible Gelfand-Tsetlin module $M$ is said to be {\em
extended\/} from $\bm\in \Specm \Gamma$ if $M(\bm)\neq 0$. The set
of isomorphism classes of irreducible Gelfand-Tsetlin modules
extended from $\bm$ is called the \emph{fiber} of $\bm\in \Specm
\Gamma$. Equivalently, this is the set of left maximal ideals of
$W(\pi)$ containing $\bm$. An important problem in the theory of
Gelfand-Tsetlin modules is to determine the cardinality of the
fiber of an arbitrary $\bm$. In the case where the fibers consist
of single isomorphism classes, the corresponding irreducible
Gelfand-Tsetlin modules are parameterized by the elements of $\Sp
\Gamma$. This problem was solved in the particular cases of
one-column pyramids \cite{O} ($\gl_n$ case) and two-row
rectangular pyramids \cite{fmo} (Yangian for $\gl_2$). We extend
these results to arbitrary finite $W$-algebras of type $A$. The
technique used in this paper is quite different, it is based on
the properties of the {\it Galois orders\/} developed in the
papers \cite{fo-Ga1} and \cite{fo-Ga2}. Our second main result is
the following theorem.

\medskip

\noindent{\bf Theorem II.} \textit{The fiber of  any $\bm\in
\Specm \Ga$ in the category of Gelfand-Tsetlin modules over
$W(\pi)$ is non-empty and finite.}

\medskip

Clearly, the same irreducible Gelfand-Tsetlin module can be
extended from different maximal ideals of $\Ga$; such ideals are
called \emph{equivalent}. Hence, Theorem~II provides a
parametrization of finite families of irreducible Gelfand-Tsetlin
modules over $W(\pi)$ by the equivalence classes of characters of
the Gelfand-Tsetlin subalgebra. Moreover, this gives a
classification of the irreducible {\it generic} Gelfand-Tsetlin
modules. In order to formulate the result, recall that a non-empty
set $X\subset \Specm \Gamma$ is called {\it massive\/} if $X$
contains the intersection of countably many dense open subsets. If
the field ${\Bbbk}$ is uncountable, then a massive set $X$ is
dense in $\Specm \Gamma$.

\medskip

\noindent{\bf Theorem III.} \textit{There exists a massive subset
$\Omega\subset \Specm \Ga$ such that} \begin{itemize}
\item[(i)]\label{item-generic-unique-module}
 \textit{For any
$\bm\in \Omega$,  there exists a unique, up to isomorphism,
irreducible module $L_{\bm}$ over $W(\pi)$ in the fiber of $\bm$}.
\item[(ii)]\label{item-generic-indecomp} \textit{For any $\bm\in
\Omega$ the extension category generated by $L_{\bm}$ contains all
indecomposable modules whose support contains $\bm$ and is
equivalent to the category of modules over the algebra of formal
power series in $n\ts p_1+(n-1)\ts p_2+\ldots+ p_n$ variables}.
\end{itemize}

We also make the following conjecture about the size of fiber in
general:

\medskip

\noindent{\bf Conjecture 1.}
 Let $(p_1, \ldots, p_n)$
 be the rows of $\pi$. For any $\bm\in \Specm \Ga$ the fiber of
$\bm$ consists of at most $p_1!(p_1+p_2)!\ldots (p_1+\ldots+
p_{n-1})!$ isomorphism classes of irreducible Gelfand-Tsetlin
$W(\pi)$-modules. The same bound holds for the dimension of the
subspace of $\bm$-nilpotents $V(\bm)$ in any irreducible
Gelfand-Tsetlin module $V$.

\medskip

This conjecture follows immediately from Theorem~5.3(iii) in
\cite{fo-Ga2} and the following conjecture.

\medskip

\noindent{\bf Conjecture 2.} $W(\pi)$ is free as left (right)
module over the Gelfand-Tsetlin subalgebra.

\medskip

These conjectures known to be true in the particular cases of
one-column pyramids \cite{O}  and two-row rectangular pyramids
\cite{fmo}.  We prove both conjectures for arbitrary two-row
pyramids (i.e., finite $W$-algebras associated with $\gl_2$).

\section{Shifted Yangians, finite $W$-algebras
and their representations}\label{sec:sy}\setcounter{equation}{0}

As in \cite{BK1}, given a pyramid
$\pi$ with the rows
$p_1\leqslant\dots\leqslant p_n$, introduce
the corresponding {\it shifted Yangian\/} $\Y_{\pi}(\gl_n)$
as the associative algebra over ${\Bbbk}$ defined by generators
\begin{align}\label{gener}
d^{\ts(r)}_i,&\quad i=1,\dots,n,&&\quad r\geqslant 1,\\
f^{(r)}_i,&\quad i=1,\dots,n-1,&&\quad r\geqslant 1,
\non\\
e^{(r)}_i,&\quad i=1,\dots,n-1,&&\quad r\geqslant p_{i+1}-p_i+1,
\non
\end{align}
subject to the following relations:
\begin{align}
[d_{i}^{\ts(r)},d_{j}^{\ts(s)}]&=0,
\non\\
[e_{i}^{(r)},f_{j}^{(s)}]&=-\ts \de_{ij}\ts\sum_{t=0}^{r+s-1}
d_{i}^{\tss\prime\ts(t)}\ts d_{i+1}^{\ts(r+s-t-1)},\non\\
[d_{i}^{\ts(r)},e_{j}^{(s)}]&=(\de_{ij}-\de_{i,j+1})\ts\sum_{t=0}^{r-1}
d_{i}^{\ts(t)}\ts e_{j}^{(r+s-t-1)},\non\\
[d_{i}^{\ts(r)},f_{j}^{(s)}]&=
(\de_{i,j+1}-\de_{ij})\ts\sum_{t=0}^{r-1}
f_{j}^{(r+s-t-1)}\ts d_{i}^{\ts(t)},
\non
\end{align}
\begin{align}
[e_{i}^{(r)},e_{i}^{(s+1)}]-[e_{i}^{(r+1)},e_{i}^{(s)}]&=
e_{i}^{(r)}e_{i}^{(s)}+e_{i}^{(s)}e_{i}^{(r)},\non\\
[f_{i}^{(r+1)},f_{i}^{(s)}]-[f_{i}^{(r)},f_{i}^{(s+1)}]&=
f_{i}^{(r)} f_{i}^{(s)}+f_{i}^{(s)} f_{i}^{(r)},\non\\
[e_{i}^{(r)},e_{i+1}^{(s+1)}]-[e_{i}^{(r+1)},e_{i+1}^{(s)}]&=
-e_{i}^{(r)}e_{i+1}^{(s)},\non\\
[f_{i}^{(r+1)},f_{i+1}^{(s)}]-[f_{i}^{(r)},f_{i+1}^{(s+1)}]&=
-f_{i+1}^{(s)}f_{i}^{(r)},
\non
\end{align}
\begin{alignat}{2}
[e_{i}^{(r)},e_{j}^{(s)}]&=0\qquad&&\text{if}\quad |i-j|>1,\non\\
[f_{i}^{(r)},f_{j}^{(s)}]&=0\qquad&&\text{if}\quad |i-j|>1,\non\\
[e_{i}^{(r)},[e_{i}^{(s)},e_{j}^{(t)}]]
&+[e_{i}^{(s)},[e_{i}^{(r)},e_{j}^{(t)}]]=0
\qquad&&\text{if}\quad |i-j|=1,\non\\
[f_{i}^{(r)},[f_{i}^{(s)},f_{j}^{(t)}]]
&+[f_{i}^{(s)},[f_{i}^{(r)},f_{j}^{(t)}]]=0
\qquad&&\text{if}\quad |i-j|=1,
\non
\end{alignat}
for all admissible $i,j,r,s,t$, where $d_{i}^{(0)}=1$ and
the elements
$d_{i}^{\tss\prime\ts(r)}$ are found from the relations
\ben
\sum_{t=0}^r d_{i}^{\tss(t)}\ts d_{i}^{\tss\prime\ts(r-t)}=\de_{r0},
\qquad r=0,1,\dots.
\een
Note that the algebra $\Y_{\pi}(\gl_n)$ depends only on
the differences $p_{i+1}-p_i$ and our definition
corresponds to the left-justified
pyramid $\pi$, as compared to \cite{BK1}.
In the particular case
of a rectangular pyramid $\pi$ with $p_1=\dots=p_n$, the
algebra $\Y_{\pi}(\gl_n)$ is isomorphic to the {\it Yangian\/}
$\Y(\gl_n)$; see e.g. \cite{m:yc} for the description
of its structure and representations.
Moreover, for an arbitrary pyramid $\pi$,
the shifted Yangian $\Y_{\pi}(\gl_n)$
can be regarded as a natural subalgebra of $\Y(\gl_n)$.

Due to the main result of \cite{BK1},
the {\it finite $W$-algebra\/} $W(\pi)$,
associated to $\gl_m$ and the pyramid $\pi$, can be defined
as the quotient of $\Y_{\pi}(\gl_n)$ by the two-sided ideal
generated by all elements $d_{1}^{\tss(r)}$ with $r\geqslant p_1+1$.
We refer the reader to \cite{BK1, BK2} for the description
of the structure of the algebra $W(\pi)$,
including analogues of the Poincar\'e--Birkhoff--Witt theorem
and a construction of algebraically independent
generators of the center.

\subsection{Gelfand-Tsetlin basis for finite-dimensional
representations}

An important role in our arguments will be played
by an explicit construction of a family of
finite-dimensional irreducible representations of
$W(\pi)$, given in \cite{FMO2}. We reproduce
some of the formulas here.

Introduce formal generating series in $u^{-1}$
with coefficients in $W(\pi)$ by
\ben
\bal
d_i(u)&=1+\sum_{r=1}^{\infty} d_i^{\ts(r)}\ts u^{-r},\qquad
f_i(u)=\sum_{r=1}^{\infty} f_i^{(r)}\ts u^{-r},\\
e_i(u)&=\sum_{r=p_{i+1}-p_i+1}^{\infty} e_i^{(r)}\ts u^{-r}
\eal
\een
and set
\ben
A_i(u)=u^{p_1}\ts (u-1)^{p_2}\ts\dots (u-i+1)^{p_i}\ts
a_i(u)
\een
for $i=1,\dots,n$ with $a_i(u)=d_1(u)\ts d_2(u-1)\dots d_i(u-i+1)$,
and
\ben
\bal
B_i(u)&=u^{p_1}\ts (u-1)^{p_2}\ts
\dots (u-i+2)^{p_{i-1}}\ts (u-i+1)^{p_{i+1}}\ts a_i(u)
\ts e_i(u-i+1),\\
C_i(u)&=u^{p_1}\ts (u-1)^{p_2}\ts\dots (u-i+1)^{p_i}
\ts f_i(u-i+1)\ts a_i(u)
\eal
\een
for $i=1,\dots,n-1$. By the results of \cite{FMO2},
$A_i(u)$, $B_i(u)$, and $C_i(u)$
are polynomials in $u$, and their coefficients
are generators of $W(\pi)$. Define
the elements $a_{r}^{(k)}$ for $r=1,\dots,n$
and $k=1,\dots,p_1+\dots+p_r$ by the expansion
\ben
A_r(u)=u^{p_1+\dots+p_r}+\sum_{k=1}^{p_1+\dots+p_r}
a_{r}^{(k)}\ts u^{p_1+\dots+p_r-k}.
\een
Then the elements $a_{r}^{(k)}$ generate the Gelfand--Tsetlin
subalgebra $\Gamma$ of $W(\pi)$ defined in the Introduction.

Recall some definitions
and results from \cite{BK2} regarding
representations of $W(\pi)$. Fix an $n$-tuple
$\la(u)=\big(\la_1(u),\dots,\la_n(u)\big)$
of monic polynomials in $u$, where $\la_i(u)$
has degree $p_i$. We let $L(\la(u))$ denote
the irreducible highest weight representation
of $W(\pi)$ with the highest weight $\la(u)$.
Then $L(\la(u))$ is generated by
a nonzero vector $\xi$ (the highest vector)
such that
\begin{alignat}{2}
B_{i}(u)\ts\xi&=0 \qquad &&\text{for} \quad
i=1,\dots,n-1, \qquad \text{and}
\non\\
u^{p_i}\ts d_{i}(u)\ts\xi&=\la_i(u)\ts\xi \qquad &&\text{for}
\quad i=1,\dots,n.
\non
\end{alignat}
Write \ben \la_i(u)=(u+\la^{(1)}_i)\ts(u+\la^{(2)}_i) \dots
(u+\la^{(p_i)}_i),\qquad i=1,\dots,n. \een We will be assuming that
the parameters $\la^{(k)}_i$ satisfy the conditions: for any value
$k\in\{1,\dots,p_i\}$ we have \ben
\la^{(k)}_i-\la^{(k)}_{i+1}\in\Z_+,\qquad i=1,\dots,n-1, \een where
$\Z_+$ denotes the set of nonnegative integers. In this case the
representation $L(\la(u))$ of $W(\pi)$ is finite-dimensional. We
will only consider a certain family of representations of $W(\pi)$
by imposing the condition \ben
\la_i^{(k)}-\la_j^{(m)}\notin\Z,\qquad\text{for all}\ \ i,j
\quad\text{and all}\ \  k\ne m. \een The {\it Gelfand--Tsetlin
pattern\/} $\mu(u)$ (associated with the highest weight $\la(u)$) is
an array of rows $(\la_{r1}(u),\dots,\la_{rr}(u))$ of monic
polynomials in $u$ for $r=1,\dots,n$, where \ben
\la_{ri}(u)=(u+\la_{ri}^{(1)})\dots(u+\la_{ri}^{(p_i)}), \qquad
1\leqslant i\leqslant r\leqslant n, \een with
$\la_{ni}^{(k)}=\la_{i}^{(k)}$, so that the top row coincides with
$\la(u)$, and \ben
\la_{r+1,i}^{(k)}-\la_{ri}^{(k)}\in\Z_+\qquad\text{and}\qquad
\la_{ri}^{(k)}-\la_{r+1,i+1}^{(k)}\in\Z_+ \een for $k=1,\dots,p_i$
and $1\leqslant i\leqslant r\leqslant n-1$.

The following theorem was proved in \cite{FMO2}. It will play
a key role in the arguments below, as it allows us to realize
$W(\pi)$ as a Galois subalgebra; see sec.~\ref{fwga}.
Set
$l^{\tss(k)}_{ri}=\la^{(k)}_{ri}-i+1$.

\bth\label{thm:dgactylp} The representation $L(\la(u))$ of the
algebra $W(\pi)$ admits a basis $\{\xi_{\mu}\}$ parameterized by all
patterns $\mu(u)$ associated with $\la(u)$ such that the action of
the generators is given by the formulas \beql{aractba} A_r(u)\ts
\xi^{}_{\mu}=\la_{r1}(u)\dots\la_{rr}(u-r+1)\ts \xi^{}_{\mu}, \eeq
for $r=1,\dots,n$, and
\begin{align}\label{bractba}
B_r(-l^{(k)}_{ri})\ts \xi^{}_{\mu}&=
-\la_{r+1,1}(-l^{\tss(k)}_{ri})\dots\la_{r+1,r+1}(-l^{\tss(k)}_{ri}-r)
\ts \xi^{}_{\mu+\de_{ri}^{(k)}},\\
C_r(-l^{(k)}_{ri})\ts \xi^{}_{\mu}&=
\la_{r-1,1}(-l^{\tss(k)}_{ri})\dots\la_{r-1,r-1}(-l^{\tss(k)}_{ri}-r+2)
\ts \xi^{}_{\mu-\de_{ri}^{(k)}}, \non
\end{align}
for $r=1,\dots,n-1$, where $\xi^{}_{\mu\pm\de_{ri}^{(k)}}$
corresponds to the pattern obtained from $\mu(u)$ by replacing
$\la_{ri}^{(k)}$ by $\la_{ri}^{(k)}\pm 1$, and the vector
$\xi^{}_{\mu}$ is considered to be zero, if $\mu(u)$ is not a
pattern. \eth

Note that the action of the operators
$B_r(u)$ and $C_r(u)$ for an arbitrary value of
$u$ can be calculated by the Lagrange
interpolation formula.

\section{Skew group structure of finite $W$-algebras}\setcounter{equation}{0}

\subsection{Skew group rings}
\label{subsection-Skew-group-rings-and-skew-group-categories} Let
$R$ be a ring, $\cM$ a subgroup of $\Aut R$, and $R*\cM$ the
corresponding skew group
ring, i.e., the free left $R$-module with the basis $\cM$ and with the
multiplication
\ben
(r_{1}m_{1})\cdot (r_{2}m_{2})= (r_{1}r_{2}^{m_{1}}) (m_{1}
m_{2}),\quad m_{1},m_{2}\in \cM,\ r_{1},r_{2}\in R.
\een

If $x\in R*\cM$ and $m\in \cM$ then denote by $x_{m}$ the element of
$R$ such that $x=\sum_{m\in \cM} x_{m} m$. Set
$$\supp x=\{m\in\cM\ |\ x_{m}\ne0\}.$$
If a finite group $G$ acts by automorphisms on $R$ and by
conjugations on $\cM$ then $G$ acts on $R*\cM$.  Denote by
 $(R*\cM)^{G}$  the subring of invariants under this action.
 Then $x\in (R*\cM)^{G}$ if and only if
$x_{m^{g}}=x^{g}_{m}$ for $m\in\cM,g\in G$.

For $\vi\in \Aut R$ and $a\in R$ set $H_{\vi}=\{h\in G|\vi^{h}=\vi\}$
and
\begin{equation}
\label{equation-definition-a-phi} [a\vi]:=\sum_{g\in
G/H_{\vi}}a^{g}\vi^{g}\in (R*\cM)^{G},
\end{equation}

where the sum is taken over representatives of the cosets and does not depend on their choice.

\subsection{Galois algebras}

Let $\Ga$ be a commutative domain, $K$  the field of fractions of
$\Ga$, $K\subset L$ a finite Galois extension, $G=\Gal(L/K)$  the
corresponding Galois group, $\cM\subset \Aut L$ a subgroup. Assume
that $G$ belongs to the normalizer of $\cM$ in $\Aut L$ and $\cM\cap
G=\{e\}$.  Then $G$ acts on the skew group algebra $L*\cM$ by
authomorphisms: $(a m)^g=a^g m^g$ where the action on $\cM$ is by
conjugation. Denote by $(L*\cM)^G$  the subalgebra of $G$-invariants
in $L*\cM$.

\bde\cite{fo-Ga1} A subring
$U\subset (L*\cM)^G$ finitely generated over $\Ga$
is called a  \emph{Galois ring over $\Ga$} if
$KU=UK=(L*\cM)^G$.

\ede

We will always assume that both $\Ga$ and $U$ are $\k$-algebras
and that $\Ga$ is noetherian. In this case we will say that a
Galois ring $U$ over $\Ga$ is a \emph{Galois algebra over} $\Ga$.

Denote by $\bar{\Ga}$ the integral closure of $\Ga$ in $L$. Let
$S_{*}=\{S_{1}\subset S_{2}\subset\dots \subset S_{N}\subset
\dots\}$ be an increasing chain of finite sets. Then the growth of
$S_{*}$ is defined as
\begin{equation}
\label{equation-definition-of-growth}
\gro(S_{*})=\overline{\lim_{N\to\infty}} \log_{N}|S_{N}|.
\end{equation}
 Fix a set of
generators $\cM_{1}=\cO_{\vi_{1}}\cup$ $\dots$ $\cup
\cO_{\vi_{n}}$ of $\cM$, where $\cO_{\vi}=\{\vi^g| g\in G\}$. For
$N \geqslant 1$, let $\cM_{N}$ be the set of words $w\in\cM$ such
that $l(w)\leqslant N$, where $l$ is the length of $w$, that is
\begin{align}
\label{align-next-set} \cM_{N+1}=\cM_{N}\bigcup
\bigg(\bigcup_{\vi\in\cM_{1}}\vi\cdot\cM_{N}\bigg).
\end{align}
  Let $\cM_{*}=\{\cM_{1}\subset
\cM_{2}\subset\dots \subset \cM_{N}\subset \dots\}$.  Then the
growth of $\cM$  is by definition $\gro(\cM_{*})$, we will denote
it by $\gro(\cM)$. For a ring $R$ we will denote by $\gkdim R$ its
Gelfand-Kirillov dimension.

\bpr \label{theorem-GK-dimension-of-some-Galois}
\cite[Theorem~6.1]{fo-Ga1} Let $U\subset L*\cM$ be a Galois
algebra over noetherian $\Ga$, $\cM$ a group  of finite growth
such that for every finite dimensional $\k$-vector space $V\subset
\bar{\Ga}$ the set $\cM\cdot V$ is contained in a finite
dimensional subspace of $\bar{\Ga}$. Then
\begin{equation}
\label{equation-formula-for-GK-dimension-of-some-Galois} \gkdim
U\geqslant\gkdim \Gamma+\gro(\cM).
\end{equation}
\epr

\subsection{PBW Galois algebras}

 Let $U$ be an associative
algebra over ${\k}$, endowed with an increasing exhausting
finite-dimensional filtration $\{U_{i} \}_{i\in\mathbb Z}$,
$U_{-1} =$ $\{ 0 \}$, $U_{0} =$ ${\k}$. Then $U_iU_j\subset
U_{i+j}$ and $\gr\ts U=$ $\displaystyle \bigoplus_{i=0}^{\infty}
U_{i }/U_{i-1 }$ is the associated graded algebra.  An algebra $U$
is called a \emph{PBW algebra} if  $\gr\ts U$ is a commutative
affine $\k$-algebra. In particular, $U$ is a noetherian affine
$\k$-algebra. For PBW algebras we have the following sufficient
conditions to be a Galois algebra.

\bth\label{thm-GKdim-PBW algebras}\cite[Theorem 7.1]{fo-Ga1} Let
$U$ be a PBW  algebra generated by the elements $u_1, \ldots, u_k$
over $\Ga$, $\gr\ts U$  a polynomial ring in $n$ variables,
$\cM\subset \Aut L$ a group and $f:U\rightarrow (L*\cM)^G$ a
homomorphism such that $\underset{i}\cup\supp f(u_i)$ generates
$\cM$. If $$\gkdim \Ga+ \gro(\cM)=n$$ then $f$ is an embedding and
$U$ is a Galois algebra over $\Ga$. \eth

\subsection{Finite $W$-algebras as  Galois algebras}
\label{fwga}

Now consider the Gelfand--Tsetlin subalgebra $\Gamma$
of the algebra $W(\pi)$, as defined in sec.~\ref{sec:sy}.
Let $\Lambda$ be the polynomial algebra in the variables
$x_{ri}^{k}, 1\leqslant i\leqslant r\leqslant n$, $k=1, \ldots,
p_i$. Consider the $\Bbbk$-homomorphism $\imath: \Gamma\to \Lambda$
defined by \beql{imath}
\imath(a_{r}^{(k)})=\sigma_{r,k}(x_{r1}^{1},\ldots, x_{r1}^{p_1},
\ldots, x_{rr}^{1},\ldots, x_{rr}^{p_r}), \qquad
k=1,\dots,p_1+\dots+p_r,
\end{equation}
where $\sigma_{r,j}$ is the $j$-th elementary symmetric polynomial
in $p_1+\ldots + p_r$ variables. The map $\imath$ is  injective by
the theory of symmetric polynomials, and we will identify the
elements of $\Gamma$ with their images in $\Lambda$. Let
$G=S_{p_1}\times S_{p_1+p_2}\times \ldots \times S_{p_1+\ldots+
p_n}$. Then $\Ga$ consists of the invariants in $\Lambda$ with
respect to the natural action of $G$. Set $\mL=\Sp \Lambda$ and
identify it with $\Bbbk^s$, $s=np_1+(n-1)p_2+\ldots +p_n$.

Let $\cM\subseteq \mL$, $\cM\simeq \Z^{(n-1)p_1+\ldots+p_{n-1}}$, be
the free abelian group generated by the symbols $\delta_{ri}^{k}\in
\Bbbk^{(n-1)p_1+\ldots+p_{n-1}}$ for $k=1, \ldots, p_i$, $1\leqslant
i\leqslant r\leqslant n-1$. Define an action
of $\cM$  on $\mL$ by the shifts
$\delta_{ri}^{k}(\ell):=\ell+\delta_{ri}^{k}$ so that
$x_{ri}^{k}$ is replaced with $x_{ri}^{k}+1$, while all other
coordinates remain unchanged. The group $G$ acts on $\mL$ by permutations and on $\cM$ by conjugations.

Let $K$ be the field of fractions of $\Ga$, $L$ the field of
fractions of $\Lambda$. Then $K\subset L$ is a finite Galois
extension with the Galois group $G$, $K=L^G$. Similarly as above
one defines the action of $\cM$ on $L$, the skew group algebra
$L*\cM$ and its invariant subalgebra $(L*\cM)^G$.

Recall the polynomials $A_i(u)$, $B_k(u)$, $C_k(u)$ in $u$,
$i=1,\dots,n$ and $k=1,\ldots, n-1$,
with coefficients in $W(\pi)$ which were
defined in Section 2.1.
Consider the polynomials $\tilde A_i(u)$, $\tilde B_k(u)$, $\tilde
C_k(u)$ in $u$, which are obtained by replacing the nonzero
coefficients of the polynomials $A_i(u)$, $B_k(u)$, $C_k(u)$
by independent variables in such a way that
the new polynomials have the same
degrees as their respective counterparts and
the polynomials $\tilde A_i(u)$ are monic.
Introduce free associative algebra T
over $\k$ generated by the coefficients of these polynomials.
 Let
$L[u]*\cM$ be the skew group algebra over the ring of polynomials
$L[u]$ and $e$  the identity element of $\cM$. Note that $A_i(u)\in
L[u]*\cM$, $i=1,\ldots n$. Introduce an algebra homomorphism
$\mathrm t: T\longmapsto L[u]*\cM$ by the formulas

\ben
\begin{aligned} \label{equation-def-t}
\mathrm t(\wt{A}_j(u))&=A_j(u)e, \\
\mathrm t(\wt{B}_r(u))&=\sum_{(s,j)} X_{rsj}^+[u]\delta_{rj}^s, \\
\mathrm t(\wt{C}_r(u))&=\sum_{(s,j)} X_{rsj}^-[u](\delta_{rj}^s)^{-1},
\end{aligned}
\een
where
$$
X_{rsj}^{+}[u]=  -\frac {\prod_{(k,i)\neq (s,j)}(u+x_{ri}^{k})}
{\prod_{(k,i)\neq
(s,j)}(x_{ri}^k-x_{rj}^s)}\prod_{m,q}(x_{r+1,q}^m-x_{rj}^s),
$$
$$
X_{rsj}^{-}[u]=\frac {\prod_{(k,i)\neq (s,j)}(u+x_{ri}^{k})}
{\prod_{(k,i)\neq
(s,j)}(x_{ri}^k-x_{rj}^s)}\prod_{m,q}(x_{r-1,q}^m-x_{rj}^s),
$$
$j=1,\dots,r$ and $s=1,\dots,p_j$. The
products $(k,i)$ associated with the variables of the form $x_{ri}^{k}$
run over the pairs with $i=1,\dots,r$ and $k=1,\dots,p_i$.

In the following lemma we use notation
\eqref{equation-definition-a-phi}.

\ble We have
\ben
\mathrm t(\wt{B}_r(u))=[X_{r11}^{+}[u]\delta_{r1}^1],\qquad
\mathrm t(\wt{C}_r(u))=[X_{r11}^{-}[u](\delta_{r1}^1)^{-1}].
\een
In
particular, $\mathrm t$ defines a homomorphism  from $T$ to
$(L*\cM)^G$. \ele

\begin{proof}
Note that $H_{\delta_{r1}^1}\subset G$ consists of permutations of
$G$ which fix $1$, and that $X_{r11}^{\pm}$ are fixed points of
$H_{\delta_{r1}^1}$. Then for
 $g\in G$, such that $g(1)=p_1+\ldots +p_{i-1}+k$, $0< k\leqslant p_i$,
the equality
 $(\delta_{r1}^1)^{g}=\delta_{ri}^k$ holds
and $(X_{r11}^{\pm})^{g}=X_{rki}^{+}$, which implies the statement.
\end{proof}

Denote by $\pi:T\myto W(\pi)$ the projection defined by
$$\wt{A}_r(u)\longmapsto A_r(u), \, \,
\wt{B}_r(u)\longmapsto B_r(u), \,\,
 \wt{C}_r(u)\longmapsto C_r(u).$$

\ble\label{lemma-g-z-formulae-embeds-gl-n-in-skew-group-alg} There
exists a homomorphism of algebras $i:W(\pi)\myto (L*\cM)^G$, such
that the diagram

$$ \xymatrix{ T \ar[rr]^{\pi}\ar[dr]_{\mathrm t} &&
W(\pi)\ar[dl]^{i}\\
&(L*\cM)^G& }   $$ commutes. \ele

\begin{proof}
Let $V$ be a finite-dimensional $W(\pi)$-module with a basis
$\{\xi_{\mu}\}$. It induces a module structure over $T$ via the
homomorphism $\pi$. Moreover, due to Theorem~\ref{thm:dgactylp}, $V$
has a right module structure over  $\mathrm t(T)\subset (L*\cM)^G$.
If $z\in T$ and $\mathrm t(z)=\ds\sum_{i=1}^{s}[a_{i}m_{i}]$,
$m_{i}\in \cM$, $a_{i}\in L$, then $ \xi_{\mu}\cdot \mathrm
t(z)=\ds\sum_{i=1}^{s}a_{i}(\mu) \xi_{m_{i}+\mu}$, where $a_i(\mu)$
means the evaluation of the rational function $a_i\in L$ in $\mu$.
Suppose now that $z\in \Ker \pi$ and consider $\mathrm t(z)$.
 There exists a dense subset
$\Omega(z)$ consisting of $\mu$'s, such that $\xi_{\mu}$ is a basis
vector of some finite-dimensional $W(\pi)$-module $V$ and
$\xi_{\mu}\cdot\mathrm t(z)$ is defined. Moreover, for any $\mu\in
\Omega(z)$, $\xi_{\mu}\cdot\mathrm t(z)=0$ and hence $a_{i}(\mu)=0$
for all $i$. Since each $a_{i}$ is a rational function on $\Specm
\Lambda$, it implies that $a_{i}=0$, and hence $z\in\Ker \mathrm t$.
Therefore, there exists a homomorphism $i:W(\pi)\myto (L*\cM)^G$
such that the diagram commutes.

\end{proof}

\bth\label{yangian-Galois} $W(\pi)$ is a Galois algebra over
$\Ga$. \eth
\begin{proof}
First note that $W(\pi)$ is a PBW algebra and $\dim_{\k} \cM\cdot
v<\infty$ for any $v\in \Lambda$. Also,
$$\gkdim W(\pi)=(2n-1)p_1+(2n-3)p_2+\ldots+3p_{n-1}+p_n=$$
$$=\gkdim \Ga +
\gro{\cM}.$$ Since $\cup_r \supp \mathrm{t}(\wt{B}_r(u))$ and
$\cup_r \supp \mathrm{t}(\wt{C}_r(u))$ contain all the generators of
the group $\cM$, all conditions of Theorem~\ref{thm-GKdim-PBW
algebras} are satisfied. Hence we conclude that $i:W(\pi)\myto
(L*\cM)^G$ is an embedding and $W(\pi)$ is a Galois algebra over $\Ga$.
\end{proof}

Recall that a commutative subalgebra $A$ of an associative algebra
$B$ is called a \emph{Harish-Chandra subalgebra} if for any $b\in
B$, the $A$-bimodule $A b A$ if finitely generated both as a left
and as a right $A$-module \cite{dfo:hc}.
\bco\label{corol-Ga-Harish-Chandra} $\Ga$ is a Harish-Chandra
subalgebra of $W(\pi)$. \eco

\begin{proof}
Since $\cM\cdot \Lambda\subset \Lambda$ and $W(\pi)$ is a Galois
algebra over $\Ga$, the statement follows from
\cite[Proposition~5.2]{fo-Ga1}.
\end{proof}

Let $\imath: K\rightarrow L$ be a canonical embedding, $\phi\in \Aut
L$, $\jmath=\phi \imath$. Consider a $(K,L)$-bimodule
$\wt{V}_{\phi}=KvL$, where $av=v\phi(a)$ for all $a\in K$. Let
 $V_{\phi}$ be the set of $\St(\jmath)$-invariant elements of
 $\wt{V}_{\phi}$.

\bco\label{corol-Galois-bimodule} Let $S=\Ga\setminus \{0\}$. Then
\begin{itemize}
\item[(i)]  $S$ is an Ore set and
$$W(\pi)[S^{-1}]\simeq (L*\cM)^G\simeq [S^{-1}]W(\pi).$$
\item[(ii)]  $K\otimes_{\Ga}W(\pi)
\otimes_{\Ga}K\simeq (L*\cM)^G$ as $K$-bimodules.
\item[(iii)]
 $W(\pi)[S^{-1}]\simeq
 \bigoplus_{\phi\in:\cM/G}V_{\phi}$ as $K$-bimodules.
\end{itemize}
\eco

\begin{proof}
Follow from Theorem~\ref{yangian-Galois} and \cite[Theorem
3.2\tss(5)]{fo-Ga1}.

\end{proof}

\section{Noncommutative Noether problem}\setcounter{equation}{0}
If $A$ is a noncommutative domain that satisfies the Ore
conditions then it admits the skew field of fractions which we
denote $D(A)$.

The $n$-th Weyl algebra  $A_{n}$ is generated by
$x_{i},\partial_{i}$, $i=1,\dots, n$ subject to relations
\begin{align}
\label{equation-commutation-rules} \nonumber
x_{i}x_{j}=x_{j}x_{i},\\
\partial_{i}\partial_{j}
=\partial_{j}\partial_{i},\\
\partial_{i}x_{j}-x_{j}\partial_{i}=\delta_{ij},\ i,j=1,\dots, n.
\end{align}

This algebra is a simple noetherian domain with the skew field of
fractions $D_n=D(A_n)$. The symmetric group $S_n$ acts on $A_n$
and hence on $D_n$ by simultaneous permutations of $x_i$'s and
$\partial_i$'s.

In this section we prove the noncommutative Noether problem for
$S_n$:

\bth\label{theorem-invariants-Weyl} $D_{n}^{S_{n}}\simeq D_{n}.$
\eth

\subsection{Symmetric differential operators}
 If
$P=\Bbbk[x_{1},\dots,x_{n}]$ then we  identify  the Weyl algebra
$A_{n}$ with the ring of differential operators $\mathcal D(P)$ on
$P$ by identifying $x_{i}$ with the operator of multiplication on
$x_{i}$ and $\partial_{i}$ with the operator of partial derivation
by $x_{i}$, $i=1,\dots, n$. If $A$ is a localization of $P$ then
$\mathcal D(A)$ is generated over $A$ by $\pa_{1},\dots,\pa_{n}$
subject to obvious relations.

It is well known that $A_n^{S_n}$ is not isomorphic $A_n$ and
hence $\mathcal D(P)^{S_n}$ is not isomorphic to $\mathcal
D(P^{S_n})$ if $n>1$. For any $i=1,\dots, n$ let  $\sigma_{i}$
denote the $i$-th symmetric polynomial in the variables
$x_{1},\dots,x_{n}$. Then
$P^{S_{n}}=\k[\sigma_{1},\dots,\sigma_{n}]\subset P$. Set
$\delta=\ds\prod_{1\leqslant i<j\leqslant n}(x_{i}-x_{j})$ and
$\Delta=\delta^{2}\in P^{S_{n}}$. Denote by $P_{\Delta}$ and
$P^{S_{n}}_{\Delta}$ the localizations of corresponding algebras
by the multiplicative set generated by $\Delta$. The canonical
embedding $i:P^{S_{n}}_{\Delta}\to P_{\Delta}$ induces a
homomorphism of algebras $$i_{\Delta}:\mathcal
D(P_{\Delta})^{S_{n}}\to \mathcal D(P^{S_{n}}_{\Delta}).$$ Let
$\mathbb A^{n}$ be the $n$-dimensional affine space over $\k$. The
algebra $\mathcal D(P_{\Delta})$ is just the ring of differential
operators on $X=\Specm P_{\Delta}\subset \mathbb A^{n}$ which is
open and $S_{n}$-invariant. The geometric quotient $X/S_{n}=\Specm
\Bbbk[\sigma_{1},\dots,\sigma_{n}]_{\Delta}$ is rational and the
projection $X\rightarrow X/S_{n}$ is etale. Since the action of
$S_n$  on $X$ is free,  $i_{\Delta}$ is an isomorphism.

\bpr The following isomorphisms hold \begin{itemize}
\item[(i)]\label{enumerate-facts-about-diff-op-local} If $A$ is a
domain, $S\subset A$ is an Ore subset then  $D(A_{S})\simeq D(A)$.

\item[(ii)]\label{enumerate-facts-about-diff-op-loc-by-discriminant}
$\mathcal D(P_{\Delta})^{{S_{n}}}\simeq (\mathcal
D(P)^{S_{n}})_{\Delta}$.

\item[(iii)]\label{enumerate-facts-about-diff-op-loc-by-dis-commuting}
$(P^{S_{n}})_{\Delta}\simeq (P_{\Delta})^{S_{n}}$.

\item[(iv)]\label{enumerate-facts-about-diff-op-loc-by-dis-comm-on-polinoms}
$\mathcal D(P_{\Delta})^{S_{n}}\simeq \mathcal
D(P^{S_{n}})_{\Delta}$.
\end{itemize}
\epr

\begin{proof}
The first statement is obvious. Note that $\mathcal D(P_S)\simeq
\mathcal D(P)_S$ for a multiplicative set $S$, \cite[Theorem
15.1.25]{mcr:nnr}. If $d\in \mathcal D(P_{\Delta})^{{S_{n}}}$ then
$d_1=\Delta^{k}d\in \mathcal D(P)^{S_n}$ for some $k\geqslant 0$
implying (ii). The third statement is obvious and (iv) follows
from the previous statements.
\end{proof}

\subsection{Proof of Theorem~\ref{theorem-invariants-Weyl}}

\begin{align}\nonumber
\label{align-why} D^{S_{n}}_{n}\simeq D^{S_{n}} (\mathcal D(P))
\simeq D^{S_{n}}(\mathcal D(P)_{\Delta})\simeq D^{S_{n}}(\mathcal
D(P_{\Delta}))  \simeq D(\mathcal D^{S_{n}}(P_{\Delta}))\simeq
D(\mathcal D((P_{\Delta})^{S_{n}}))\\
\nonumber \simeq D(\mathcal
D(\k[\sigma_{1},\dots,\sigma_{n}]_{\Delta})) \simeq D((\mathcal
D(\k[\sigma_{1},\dots,\sigma_{n}])_{\Delta})\simeq D(\mathcal
D(\k[\sigma_{1},\dots,\sigma_{n}]))\simeq D_{n}.
\end{align}

 Hence $D_{n}^{S_{n}}\simeq
D_{n}$.

\section{Gelfand-Kirillov conjecture}
\setcounter{equation}{0}

Since $W(\pi)$ is a noetherian integral domain with a polynomial
graded algebra, then it satisfies the Ore conditions by the Goldie
theorem. Let $D_{\pi}(n)=D(W(\pi))$ be the skew field of fractions
 of $W(\pi)$. Recall the structure of $W(\pi)$ as a Galois algebra
over $\Ga$: $W(\pi)\subset (L*\cM)^G$, where $L$ is a field of
rational functions in $x_{ij}^k$, $j=1, \ldots, i, k=1, \ldots,
p_i$, $i=1, \ldots, n$. Then $D_{\pi}(n)\simeq D((L*\cM)^G)$.
Moreover, we will see below that $L*\cM$ has a skew field of
fractions and thus $D_{\pi}(n)\simeq D(L*\cM)^G$ \cite[Theorem
1]{fa}. Since $\Ga$ is a Harish-Chandra subalgebra
(Corollary~\ref{corol-Ga-Harish-Chandra}) then by
 \cite[Theorem 4.1]{fo-Ga1},  we have

\bpr\label{prop-skew-center} The center $\cZ$ of $D_{\pi}(n)$ is
isomorphic to $K^{\cM}$. \epr

Let $\Lambda$ be the polynomial ring in variables
$x_{ij}^k$, $j=1, \ldots, i, k=1, \ldots, p_j$, $i=1, \ldots, n$.
Denote  by $L_i$ (respectively $\Lambda_i$) the field of
rational functions (respectively the polynomial ring)
in $x_{ij}^k$ with fixed $i$. Then
$$\Lambda*\cM^G\simeq \otimes_{i=1}^{n-1}
(\Lambda_{i}*\mathbb Z^{p_1+\ldots+p_i})^{S_{p_1+\ldots+p_i}}
\otimes \Lambda_n^{S_{p_1+\ldots+p_n}}.
$$

\bpr For every $i=1, \ldots, n$
$$D(L_i*\mathbb Z^{p_1+\ldots+p_i})
\simeq D(A_{p_1+\ldots+p_i}(\Bbbk)).$$
\epr

\begin{proof}
Let $B_i=\Bbbk[t_1, \ldots, t_i]*\mathbb Z^i$, where $\mathbb Z^i$
is generated by $\sigma_k$, $k=1, \ldots, i$ and
$\sigma_k(t_m)=t_m-\delta_{km}$. Then $B_i$ is isomorphic to
 the localization $\cA_i$ of the $i$-th
 Weyl algebra with respect to $x_1, \ldots, x_i$.
 This isomorphism is given as follows:
$$x_k\mapsto \sigma_k, \,\,\, \partial_k\mapsto t_k\sigma_k^{-1}.$$
Hence, a subring $\Lambda_i*\mathbb Z^{p_1+\ldots+p_i}$ of
$L_i*\mathbb Z^{p_1+\ldots+p_i}$  is isomorphic to a localization
of $A_{p_1+\ldots+p_i}(\Bbbk)$ which implies the statement.
\end{proof}

Since $D(A_k)^{S_k}\simeq D(A_k^{S_k})$  then we have the
isomorphism
$$D((L*\cM)^G)=D((\Lambda*\cM)^G)\simeq
\otimes_{i=1}^{n-1}D((A_{p_1+\ldots+p_i}
(\Bbbk))^{S_{p_1+\ldots+p_i}}\otimes D(T_n)),
$$
where $T_n=\Lambda_n^{S_{p_1+\ldots+p_n}}$ is a polynomial ring
isomorphic $\Lambda_n$. Moreover, by applying
Theorem~\ref{theorem-invariants-Weyl} we have the isomorphism
$$D((L*\cM)^G) \simeq D(A_{(n-1)p_1+\ldots +p_{n-1}}(\Bbbk)\otimes D(T_n)).$$

Since $D(T_n)$ is a pure transcendental extension of $\Bbbk$ of
degree $p_1+\ldots+p_n$, and since $D((L*\cM)^G)\simeq D(W(\pi))$,
 we have thus proved the Gelfand-Kirillov conjecture (Theorem I):
$$D(W(\pi))\simeq D(A_{(n-1)p_1+\ldots +p_{n-1}}
(D(T_n)))=D_{k, m},$$
$k=(n-1)p_1+\ldots +p_{n-1}$, $m=p_1+\ldots+p_n$.

Recall that the {\it Miura transform\/} \cite{BK2} is an injective
homomorphism
$$\tau:W(\pi)\rightarrow \otimes_{i=1}^l U(\gl_{q_i}).$$
Observe that
$D(\otimes_{i=1}^l U(\gl_{q_i}))\simeq D_{k,m}$, since
$k=\sum_{i=1}^l q_i(q_i-1)/2$ and $m=\sum_{i=1}^l q_i$.
Hence we have proved the following corollary.

\bco\label{corol-Miura-transform} The Miura transform extends to an
isomorphism of the corresponding skew fields of fractions.
\eco

\section{Fibers of characters}
\setcounter{equation}{0}

\subsection{Galois orders}

Let $U\subset (L*\cM)^G$ be a Galois ring over an integral domain
$\Gamma$.

 \bde\cite{fo-Ga1} A Galois ring $U$  over   $\Ga$ is called
a \emph{Galois order}{}  if for any finite dimensional right
(respectively left) $K$-subspace $V\subset U[S^{-1}]$
(respectively $V\subset [S^{-1}]U$), $V\cap U$ is a finitely
generated right (respectively left) $\Ga$-module. \ede A concept
of a Galois order over $\Ga$ is a natural noncommutative
generalization of a classical notion of $\Ga$-order in skew group
ring $(L*\cM)^G$. If $\Ga$ is a noetherian $\k$-algebra then a
Galois order over $\Ga$ will be called an \emph{integral Galois
algebra}. Note that in particular a Galois ring $U$ over $\Ga$ is
a Galois order if $U$ is a projective right and left $\Ga$-module.

The following criterion  for Galois orders was established in
\cite[Corollary 5.6]{fo-Ga1}.

\bpr \label{prop-criterion-integra-algebra} Let $U\subset L*\cM$ be
a Galois algebra over a noetherian normal $\Bbbk$-algebra $\Ga$.
Then the following statements are equivalent
\begin{itemize}
\item[(i)]
\label{enum-algebra-is-integral} $U$ is an integral Galois algebra over
$\Ga$. \item[(ii)] \label{enum-divisor-criterion-integrality} $\Ga$
is a Harish-Chandra subalgebra and, if for $u\in U$ there exists a
nonzero $\ga\in \Ga$ such that $\ga u\in \Ga$ or $u\ga\in \Ga$, then
$u\in\Gamma$.
\end{itemize}
\epr

Suppose now that $U$ is a PBW Galois algebra over $\Ga$
 with the polynomial associated graded algebra
$\gr\ts U=A$. Then both $U$ and $A$ are endowed with  degree
function $\deg$ with obvious properties. For $u\in U$ denote by
$\bar{u}\in A$ the corresponding homogeneous element. Also denote by
$\gr \Ga$ the image of $\Ga$ in $A$. Then we have the following
graded version of Proposition~\ref{prop-criterion-integra-algebra}.

\ble \label{lemma-reduction-to-commutative} Let $U\subset L*\cM$ be
a PBW Galois algebra over a noetherian normal $\Bbbk$-algebra $\Ga$
with a polynomial graded algebra $\gr\ts U$. Then the following
statements are equivalent
\begin{itemize}
\item[(i)]
\label{enum-algebra-is-integral-graded} $U$ is an integral Galois
algebra over $\Ga$.
\item[(ii)]\label{enumerate-just-reduction-to-commutative}
$\Ga$ is a Harish-Chandra
subalgebra and for $\gamma,\gamma'\in \Ga\setminus \{0\}$
it follows from
$\bar{\ga'}=\bar{\ga} a, a\in A$ that $a\in \gr \Ga$.
\end{itemize}
\ele

\begin{proof}
Suppose $\ga'=\ga u\ne0$, $\ga',\ga\in \Ga$, $u\in U\setminus \Ga$
and $\deg \ga'$ is the minimal possible. Then $\bar{\ga'}=\bar{\ga}
\bar{u}\ne0$ in $A$. By the assumption $\bar{u}=\bar{\ga}''$ for
some in $\ga''\in \Ga$ and hence either $\ga''=u$, or $\ga_{2}=\ga
u_{1}\in \Ga$, where $u_{1}=u-\ga''$, $\ga_{2}=\ga'-\ga \ga''$.
Since in the second case $\deg \ga_{2}<\deg \ga_{1}$ this
contradicts the minimality assumption. Therefore, $\ga''=u\in \Ga$.
The case  $\ga'=u\ga \ne0$ is considered analogously. Hence the
statement \eqref{enum-divisor-criterion-integrality} of
Proposition~\ref{prop-criterion-integra-algebra} holds, which
implies the integrality of the Galois algebra $U$.
 \end{proof}

Representation theory of Galois algebras was developed in
\cite{fo-Ga2}. For $\bm\in \Specm \Ga$ denote by $F(\bm)$ the
fiber of $\bm$ consisting of isomorphism classes of irreducible
Gelfand-Tsetlin (with respect to $\Ga$) $U$-modules $M$ with
$M(\bf m)\neq 0$.

Let $E$ be the integral extension of $\Ga$ such that $\Ga=E^G$ and
assume that $\Ga$ is noetherian. Then the fibers of the surjective
map $\vi:\Specm E\rightarrow \Specm \Ga$ are finite. Let $\bm\in
\Specm \Ga$ and $l_{\bm}\in \Specm E$ such that $\vi(l_{\bm})=\bm$.
Denote
$$\St_{\cM}(\bm)=\{x\in \cM|x\cdot l_{\bm}=l_{\bm}\}.$$
Clearly the set $\St_{\cM}(\bm)$ does not depend on the choice of
$l_{\bm}$. \bth\label{thm-Galois-representations}\cite[Theorem
A]{fo-Ga2} Let $U$ be an integral Galois algebra over noetherian
$\Ga$, $\bm\in \Specm \Ga$. If the set $\St_{\cM}(\bm)$ is finite
then the fiber $F(\bm)$ is non-trivial and finite. \eth

\subsection{Finite $W$-algebras as integral Galois algebras}

In this section we show that $W(\pi)$ is an integral Galois
algebra over $\Ga$.

Following \cite[Section~2.2]{BK2}, for
$1\leqslant i\leqslant j\leqslant n$
define the higher root elements $e_{ij}^{(r)}$
and $f_{ji}^{(r)}$ of $W(\pi)$ inductively by the formulas
$e_{i,i+1}^{(r)}=e_i^{(r)}$ for $r\geqslant p_{i+1}-p_i+1$,
\ben
e_{ij}^{(r)}
=[e_{i,j-1}^{(r-p_j+p_{j-1})},e_{j-1}^{(p_j-p_{j-1}+1)}]
\qquad\text{for}\quad r\geqslant p_j-p_i+1,
\een
and
\ben
f_{i+1,i}^{(r)}=f_i^{(r)},\qquad f_{j,i}^{(r)}
=[f_{j-1}^{(1)},f_{j-1,i}^{(r)}]
\qquad\text{for}\quad r\geqslant 1.
\een
Furthermore, set
\ben
e_{ij}(u)=\sum_{r=p_j-p_i+1}^{\infty}e_{ij}^{(r)}\ts u^{-r},
\qquad
f_{ji}(u)=\sum_{r=1}^{\infty}f_{ji}^{(r)}\ts u^{-r},
\een
and define a power series
\ben
t_{ij}(u)=\sum_{r\geqslant 0} t_{ij}^{(r)}\ts u^{-r}
=\sum_{k=1}^{\min\{i,j\}} f_{ik}(u)\ts d_k(u)\ts e_{kj}(u)
\een
for some elements $t_{ij}^{(r)}\in W(\pi)$. Due to
\cite[Lemma~3.6]{BK2}, an ascending filtration
on $W(\pi)$ can be defined by setting $\deg t_{ij}^{(k)}=k$.
Let $\overline{W}(\pi)=\gr\ts W(\pi)$ denote
the associated graded algebra and let
$\overline{t}_{ij}^{\ts(r)}$ denote the image of $t_{ij}^{(r)}$
in the $r$th component of $\gr\ts W(\pi)$. Then
$\overline{W}(\pi)$ is a polynomial algebra in
the variables
\ben
\overline {t}_{ij}^{\ts(r)}\quad\text{with}\quad i\geqslant j,\quad
1\leqslant r\leqslant p_j\quad\text{and}\quad
\overline {t}_{ij}^{\ts(r)}\quad\text{with}\quad i< j,\quad
p_j-p_i+1\leqslant r\leqslant p_j.
\een
By \cite[Theorem~3.5]{BK2}, the series
\ben
T_{ij}(u)=u^{p_j}\ts t_{ij}(u),\qquad 1\leqslant i,j\leqslant n,
\een
are polynomials in $u$. Introduce
the matrix
$T(u)=(T_{ij}(u-j+1))_{i,j=1}^n$ and consider its
\emph{column determinant}
\begin{equation}
\begin{split}
\cdet T(u)= \sum_{\sigma\in S_n} \sgn\ts\sigma\cdot
T_{{\sigma(1)}1}(u)T_{{\sigma(2)}2}(u-1)...T_{{\sigma(n)}n}(u-n+1).
\end{split}
\end{equation}
This is a polynomial in $u$, and the  coefficients ${d}_{s}\in
W(\pi)$ of the powers $u^{p_1+\ldots+p_n-s}$, $s=1, \ldots,
p_1+\ldots +p_n\,$ are algebraically independent generators of the
center  of $W(\pi)$; see \cite{bb:ei}.

For $F=\sum_i f_iu^i\in W(\pi)[u]$ denote $\overline F=\sum_i
\overline {f}_i u^i\in \overline {W}(\pi)[u]$. Also we denote
$X_{ij}^{k}=\overline {t}_{ij}^{(k)}$ for $k\geqslant 1$ and
$X_{ij}^0=\delta_{ij}$. Set $X_{ij}(u)=\overline {T}_{ij}(u)$ and
$X(u)=(X_{ij}(u))_{i,j=1}^n$. Since $\overline{
{T}_{ij}(u-\lambda)}=X_{ij}(u)$ for any $\lambda\in\k$, one can
easily check  that $\gr\ts\cdet T(u)=\det X(u)$.

Then
\begin{equation} \label{equations-coeff-center}
\overline{d}_{s}=\sum_{\substack{k_{1}+\dots+k_{n}=s}}\
\sum_{\sigma\in S_{n}}\sgn\ts\sigma\cdot X_{\sigma(1)1}^{k_{1}}\dots
X_{\sigma(n)n}^{k_{n}}
\end{equation}
is just the coefficient of $u^{p_1+\dots+p_n-s}$ in $\det X(u)$.

Fix $r$, $1\leqslant r\leqslant n$ and consider
$X_r(u)=(X_{ij}(u))_{i,j=1}^r$. Then
\begin{equation} \label{equations-coefficients}
{d}_{r\,s}=\sum_{\substack{k_{1}+\dots+k_{r}=s}}\ \sum_{\sigma\in
S_{r}}\sgn\ts\sigma\cdot X_{\sigma(1)1}^{k_{1}}\dots
X_{\sigma(r)r}^{k_{r}}
\end{equation}
is  the coefficient of $u^{p_1+\ldots+p_r-s}$ in $\det X_r(u)$ and
the elements
$$\{d_{r\, s}, \quad r=1, \ldots, n,\quad s=1, \ldots, p_1+ \ldots + p_r\}$$ are
the generators of the algebra $\gr\ts \Ga$.

We will use the idea of a weighted polynomial order on $\overline
{W}(\pi)$ (\cite{GP}). Let $S=\{X_{ij}^{k}\,|\,i,j=1,\dots,n;
k=1,\dots,p_j\}$, $w:S\to \bN$ be a function. Define the degree of
each variable $X_{ij}^k$ as $w(X_{ij}^k)$ and then define the
degree of any monomial in these variables as the sum of the
degrees of the variables occurring in the monomial. We will denote
this degree associated with $w$ by $\deg_{w}$.
It coincides with the usual polynomial degree if $w(X_{ij}^{k})=1$
for all $i,j,k$. Also it coincides with the degree in $\overline
{W}(\pi)$ if $w(X_{ij}^{k})=k$ for all $i,j$. Fixing an order on $S$ we define a lexicographic order on the monomials. For the monomials $m_{1}$ and $m_{2}$ define
$m_{1}>_w m_{2}$ provided  $\deg_{w}(m_{1})>\deg_{w}(m_{2})$
or $\deg_{w}(m_{1})=\deg_{w}(m_{2})$ and $m_{1}>m_{2}$ in the
lexicographical order. It allows us to define  the leading monomial of $f\in \overline {W}(\pi)$ 
with respect
to $w$.    If $m$ is a 
leading monomial of $f$ then set $\lm(f)=m$. The coefficient of
$m$ in $f$ we denote by $\lc(m)=\lc(f)$. Note that the weighted
polynomial order $\deg_{w}$ and the concepts of  $\lm(f)$
and $\lc(f)$ naturally extend to
 $W(\pi)$ and $\Gamma$.

\ble \label{lemma-existence-of-the-weight}
\begin{enumerate}

\item\label{existence-variable} There exists a weight function $w$ such that for  any $r=1,
\ldots, n$ and $s=1, \ldots ,p_1+\ldots+p_r$
the leading monomial $m_{r\,s}$ of $d_{r\,s}$  contains a variable $X(r,s)=X_{ij}^{k}$, $(i,j,k=i(r,s),j(r,s),k(r,s))$
which does not enter in leading monomials  $m_{r'\,s'}$ for $(r',s')\ne(r,s)$. Besides, the variable $X(r,s)$ enters in
 $m_{r\,s}$ in degree $1$.

 \item
\label{enumerate-monotone-groebner-defines-leading-in-gamma} For
any  $\ga\in \gr \Ga$ there exist   $f=\prod_{r,s}d_{r\, s}^{\,k_{r,s}}$ and  $\la\in \Bbbk$ such that    $\ga>_w(\ga-\la f).$

\item \label{enumerate-monotone-groebner-division-of-leading} If
$\ga,\ga_{1}\in \gr \Ga$ and $\lm(\ga_1)|\lm(\ga),$ then there
exists $\ga_{2}\in \gr \Ga,$ such that
    $$\lm(\ga)=\lm(\ga_{1})\lm(\ga_{2}).$$

\end{enumerate}

\ele

\begin{proof}
Define a function $v$ on $S$ with values in $\bZ$ satisfying the
following conditions:
\begin{itemize}
\item[(i)]\label{positive} $v(X_{i+1\,i}^{p_i})=i+1$, $i=1,\dots, n-1$;
\item[(ii)]\label{moderate}
$v(X_{ij}^{k})=-N$, where  $N> 2 n^{2}$,  if $i<j$, $i,j=1,\dots,n$;
\item[(iii)]\label{bad}
$v(X_{ii}^{k})$ are significantly smaller than those above,\\
$v(X_{ii}^{k})>v(X_{jj}^{l})$  if
 $i>j$ or  $i=j, k>l$;
\item[(iv)]\label{very-very-bad}  For $i-j\geqslant 2$ or
$j=i-1, k<p_{i-1}$ the values $v(X_{ij}^{k})$ are negative and its absolute
values are significantly larger  than the absolute
values of those above.
\end{itemize}

In particular,  if a monomial $m$ from
\eqref{equations-coefficients}  contains $X_{ij}^{k}$ satisfying (iv), then 
 $v(m)< v(m')$ for any $m'$ which does not contain such variable.
 
First we will construct a required monomial for the weight
function $v$. Fix $r\in \{1, \ldots, n\}$ and $s\in \{1, \ldots
p_1+\ldots+p_r\}$. 

If $s\leqslant p_r$ then  set
$$y_{r,s}=X_{r\,r}^s.$$ Suppose $p_r< s\leqslant p_r+p_{r-1}$ and consider
$$y_{r,s}=X_{r,r-1}^{p_{r-1}}X_{r-1,r}^{s-p_{r-1}}.$$ Note that 
$s-p_{r-1}\leqslant p_r$ and $p_r-p_{r-1}< s-p_{r-1}$. Generalizing, suppose 
$$p_{r}+\ldots +p_{r-t+1}< s\leqslant p_{r}+\ldots +p_{r-t},$$ for some $t$, 
$2\leqslant t\leqslant r-1$ (such $t$ is uniquely defined for  given $r, s$). In this case set 
$$y_{r,s}=X_{r,r-1}^{p_{r-1}}X_{r-1,r-2}^{p_{r-2}}\ldots X_{r-t+1,r-t}^{p_{r-t}}X_{r-t,r}^{k},$$
where $k=s-(p_{r-1}+\ldots +p_{r-t})$. We have $p_r-p_{r-t}< k\leqslant p_r$.

It is easy to see that the defined monomials $y_{r,s}$ belong to
${d}_{r\,s}$. Moreover, any other monomial in  ${d}_{r\,s}$ has weight strictly smaller than $y_{r,s}$.
Indeed, the condition (iv) shows that if a leading monomial
in ${d}_{r\,s}$ contains $X_{ij}^{k}$, where $i>j$, then $i=j+1$
and $k=p_j$. Hence $y_{r,s}$ is the leading monomial of
${d}_{r\,s}$ if $s\leqslant p_r$. For the case $s> p_r$ the
conditions (iii) and (iv) show that the leading monomial of
$d_{r\,s}$ contains only $X_{i+1\,i}^{p_i}$ and $X_{ij}^{b}$ for
$i<j$. By the condition (i) we have
$$v(X_{r,r-1}^{p_{r-1}})>v(X_{r-1,r-2}^{p_{r-2}})> \dots
>v(X_{2\,1}^{p_1})$$ and hence $X_{i+1\, i}^{p_i}$
will enter the leading monomial with a largest possible value of
$i$. It is clear now that $y_{r,s}$ is the leading monomial
of $d_{r\,s}$.

Now choose a sufficiently large integer $l>0$ such that
$v(x_{ij}^{k})+kl\in \bN$ for all possible $i,j,k$. We can define
the required function $w:S\rightarrow \bN$  by
$w(x_{ij}^{k})=v(x_{ij}^{k})+kl$.
 Since $d_{r\,s}$ are homogeneous, their leading monomials do
not change after the shift of gradation. We conclude that with
respect to the function $w$, the elements
$$\{y_{r,s}\,|\,r=1,\dots, n;\, s=1,\dots, p_1+\ldots +p_r\}$$ are
the leading monomials of the generators of $\gr \Ga\subset \overline{W}(\pi)$.

Note that $y_{r,s}\neq y_{r',s'}$  for different pairs $r, s$ and
$r', s'$. Given $r$ and $s$ let $t$ be such that $0\leqslant t\leqslant r-1$ and
$p_{r}+\ldots +p_{r-t+1}< s\leqslant p_{r}+\ldots +p_{r-t}$. Set
 $X(r,s)=X_{r-t, r}^k$, $k=s-(p_{r-1}+\ldots +p_{r-t})$ if $t>0$, and $X(r,s)=X_{r\,r}^s$ if $t=0$. Then $X(r,s)$ satisfies 
 \eqref{existence-variable}.

For any $\gamma \in \gr \Ga$,
 the number of occurrences of
$d_{rs}$ in $\lm(\gamma)$ equals the number of occurrences of
$X(r,s)$ in $\lm(\gamma)$. Denote this number by $k_{r\,s}$ and set 
$f=\prod_{r,s}d_{r\, s}^{\,k_{r,s}}$. Let $\lambda=\lc(\ga)$. Then 
$$\deg_w(\gamma)> \deg_w(\gamma - \lambda f),$$ 
implying
\eqref{enumerate-monotone-groebner-defines-leading-in-gamma} and
thus \eqref{enumerate-monotone-groebner-division-of-leading}.

\end{proof}

\bth \label{theorem-restricted-yangian-is-integral} Let
$\Gamma\subset W(\pi)$ be the Gelfand-Tsetlin subalgebra of
$W(\pi)$. Then $W(\pi)$ is an integral Galois algebra over $\Ga$.
\eth

\begin{proof}
First recall that $\Ga$ is a Harish-Chandra subalgebra. Assume
that $\gamma a\in \gr\ts\Gamma$ for some $\gamma\in \gr\ts  \Ga$
and $a\in \bar{W}(\pi)$. Let $w$ be the function constructed in
Lemma~\ref{lemma-existence-of-the-weight}.
 Then $\lm(\gamma a)=\lm(\gamma)\lt(a).$ Following Lemma
\ref{lemma-existence-of-the-weight},
\eqref{enumerate-monotone-groebner-division-of-leading}, there
exists $\ga'\in \gr\ts \Ga,$ such that $\lm(\ga')=\lt(a).$
Consider $a'=a-\ga'$. Then we have $\ga a'\in \gr\ts \Ga$ and
$\deg_w(a')<\deg_w(a)$. Applying induction in $\deg_w(a)$ we
conclude that $a'\in \gr\ts\Ga$ and hence $a\in \gr\ts\Ga$. It
remains to apply Lemma~\ref{lemma-reduction-to-commutative}.
\end{proof}

Since $W(\pi)$ is integral Galois algebra over $\Ga$ and
$\Ga$ is noetherian then $W(\pi)\cap K\subset L$ is an
integral extension of $\Ga$ by  \cite[Theorem 5.2]{fo-Ga1}. Since
$W(\pi)$ is a Galois algebra over $\Ga$ then $K\cap
W(\pi)$ is a maximal commutative $\k$-subalgebra in
$W(\pi)$ by \cite[Theorem~4.1]{fo-Ga1}. But $\Ga$ is
integrally closed in $K$. Hence we obtain

\bco $\Ga$ is a maximal commutative subalgebra in $W(\pi)$.
\eco

\subsection{Proof of Theorem II}
We are now in a position to prove our main result on
Gelfand-Tsetlin modules announced in Introduction. Since the
Gelfand-Tsetlin subalgebra is a polynomial ring, $W(\pi)$ is
integral Galois algebra by
Theorem~\ref{theorem-restricted-yangian-is-integral}, and  since
for any $\bm\in \Specm \Ga$ the set $\St_{\cM}(\bm)$ is finite,
then  Theorem~II follows immediately from
Theorem~\ref{thm-Galois-representations}. Therefore every
character $\chi:\Ga \rightarrow \Bbbk$ of the Gelfand-Tsetlin
 subalgebra  defines an irreducible Gelfand-Tsetlin module which is
 a quotient of $W(\pi)/W(\pi)\bm$,
$\bm = \Ker \chi$. Of course different characters can give
isomorphic irreducible modules. In such case we say that these
characters are equivalent.
  Therefore we obtain a classification of irreducible
  Gelfand-Tsetlin modules
  by the equivalence classes of characters of $\Ga$ up to a
  certain finiteness. This finiteness corresponds to  finite fibers of irreducible
  Gelfand-Tsetlin modules with a given character of $\Ga$.

\section{Category of Gelfand-Tsetlin
modules}\label{sec:cat-of-HC-mod} \setcounter{equation}{0}

For a $\Ga$-bimodule $V$ denote by $_{\mathbf n}\hat{V}_{\mathbf
m}$ the $I$-adic  completion of $\Ga\otimes_{\k}\Ga$-module $V$,
where $I\subset \Ga\otimes \Ga$ is a maximal ideal $I=\mathbf
n\otimes\Ga+\Ga\otimes \mathbf m$, that is
\begin{equation*}
 _{\mathbf n}\hat{V}_{\mathbf m} \ = \varprojlim_{n,m} {} _{\mathbf n^n}V_{
\mathbf m^m},
\end{equation*}
here ${} _{\mathbf n^n}V_{ \mathbf m^m}=V/(\mathbf n^{n}V+V\mathbf
m^{m})$.
 Let $F(W(\pi))$ be the set of finitely generated
$\Ga$-subbimodules in $W(\pi)$.

Define a category \index{${\cA}$} ${\cA}$ $=$ ${\cA}_{U,{\Gamma}}$
with the set of objects $\Ob{\cA}=$ $\Specm {\Gamma}$ and with the
space of morphisms ${\mathscr A}({\mathbf m},{\mathbf n})$ from
${\mathbf m}$ to ${\mathbf n}$, where
\begin{equation*}
 {\cA}(\mathbf m,\mathbf n) \ =
\varinjlim_{V\in F(W(\pi))} {}_{\mathbf n}\hat{V}_{\mathbf m}.
\end{equation*}

Consider the completion $\displaystyle {\Gamma}_{{\bm}}=$
$\displaystyle \lim_{\leftarrow n}{\Gamma}/{\bm}^n$ of ${\Gamma}$
by the ideal ${\bm} \in$ $\Sp {\Gamma}$. Then the space ${\mathscr
A}({\bm},{\mathbf n})$ has a natural structure of
$({\Gamma}_{{\bn}},{\Gamma}_{{\bm}})$-bimodule. The category
${\mathscr A}$ is naturally endowed with the topology of the
inverse limit. Consider the category ${\mathscr A}\dmo_d$ of
continuous functors $M:{\mathscr A}\to{\Bbbk}\dmo$,
\cite[Section~1.5]{dfo:hc}, where ${\Bbbk}\dmo$ is endowed with
the discrete topology.

Let $\mathbb H(W(\pi),{\Gamma})$ denote the category of
Gelfand-Tsetlin modules with respect to the Gelfand-Tsetlin
subalgebra $\Ga$ for finite $W$-algebra $W(\pi)$. Since $\Ga$ is a
Harish-Chandra subalgebra by
Corollary~\ref{corol-Ga-Harish-Chandra} then by \cite[Theorem
17]{dfo:hc} (see also \cite[Theorem 3.2]{fo-Ga2}), the categories
${\mathscr A}\dmo_d$ and $\mathbb H(W(\pi),{\Gamma})$ are
equivalent.

A functor that determines this equivalence can be defined as
follows. For $N\in {\mathscr A}\dmo_d$ set
\ben \mathbb
F(N)=\osu{{\bm}\in\Specm {\Gamma}}{} N({\bm})
\een
and for $x\in N({\mathbf m}),\  a\in U$ set
\ben
ax= \sum_{{\mathbf n}\in \Specm
{\Gamma}} a_{{\mathbf n}}x,
\een
where $a_{{\mathbf n}}$ is the image of $a$ in ${\mathscr
A}({\mathbf m},{\mathbf n})$. If $f:M{\longrightarrow} N$ is a
morphism in ${\mathscr A}\dmo_d$ then set $\mathbb F(f)=$
$\oplus_{{\mathbf m}\in \Specm {\Gamma}} f({{\mathbf m}})$. Hence we
obtain a functor
$$\mathbb F:{\mathscr A}\dmo_d {\longrightarrow}
\mathbb H(W(\pi),{\Gamma}).$$

 For ${\mathbf m}\in \Sp \Gamma$ denote by
$\hat{\mathbf m}$ the completion of $\mathbf m$.  Consider the
two-sided ideal $I\subseteq \mathscr A$ generated by the completions
$\hat{\mathbf m}$ for all $\mathbf m\in \Sp \Gamma$ and set
${\mathscr A}_W= {\mathscr A}/I$.

Let $\mathbb HW(W(\pi),{\Gamma})$ be the full subcategory of
\emph{weight} Gelfand-Tsetlin modules $M$ such that $\bm v=0$ for
any $v\in M(\bm)$. Clearly,
 the categories  {$ \mathbb HW(W(\pi),{\Gamma})$}{} and
${\mathscr A}_W\dmo$ are equivalent.

For a given $\bm\in \Specm \Ga$ denote by ${\mathscr A}_{\bm}$ the
indecomposable block  of the category ${\mathscr A}$ which contains
$\bm$.

An embedding $\imath: \Gamma\to \Lambda$ induces an epimorphism
$$\imath^*: \mL\rightarrow \Specm \Ga.$$ Denote by
$\wt{\Omega}\subset \mL$ the set of generic parameters
$\mu=(\mu_{ij}^k, i=1, \ldots, n; j=1, \ldots i; k=1, \ldots p)$
such that
$$\mu_{ij}^k-\mu_{i,s}^q\notin \mathbb Z, \,\,
\mu_{r+1,j}^{(m)}-\mu_{ri}^{(k)}\notin \Z$$ for all $r,i,j,m,k$.

\bth\label{thm:category-GT-modules} Let $\bm\in \Sp\Ga$,  $\mu\in
(\imath^*)^{-1}(\bm)$. Suppose $\mu\in \wt{\Omega}$. Then
\begin{itemize}
\item[(i)]
\label{enum-homs-in-ha-cha-cat-gamma-completed-in-point} All
objects of $\cA_{\bm}$ are isomorphic
 and    for every $\bn \in \cA_{\bm}$,
$$\cA(\bn,\bn)\simeq \hat{\Ga}_{\bn}.$$

\item[(ii)]\label{enum-homs-in-ha-cha-cat-one-dim-univ-mod} Let
    $M_{\bm}=\cA_{\bm}/\cA_{\bm}\hat{\bm}$.  Then there is a
    canonical isomorphism
    $$\bF(M_{\bm})\simeq W(\pi)/W(\pi){\bm}.$$

\item[(iii)]\label{enum-homs-equiv} The category $ \mathbb
H(W(\pi),{\Gamma},\bm)$ which consists of modules whose support
belongs to $\A_{\bm}$, is equivalent to the extension category
generated by module $\bF(M_{\bm})$. Moreover, this category is
equivalent to the category $\hat{\Gamma}_{\bm}\text{\rm-mod}$.
\end{itemize}
\eth

\begin{proof}
Since $\cM$ acts freely on $\wt{\Omega}$ and $\cM\cdot \mu\cap
G\cdot \mu=\{\mu\}$ all statements follow from
Theorem~\ref{theorem-restricted-yangian-is-integral} and
\cite[Theorem~5.3, Theorem~B]{fo-Ga2}.
\end{proof}

Since for $\bm$ from Theorem~\ref{thm:category-GT-modules},
$\hat{\Ga}_{\bm}$ is isomorphic to the algebra of formal power
series in $\gkdim \Ga$ variables, we immediately obtain  the
statements of Theorem~III.

\section{$W$-algebras associated with $\gl_2$}\setcounter{equation}{0}

In this section we consider the case of $W$-algebras associated
with $\gl_2$: $W(\pi)$,  where $\pi$ has rows $(p_1,p_2)$. We will
show that $W(\pi)$ is free over the Gelfand-Tsetlin subalgebra. A
particular case $p_1=p_2$ was considered in \cite{fmo}.

The shifted Yangian $W(\pi)$ is generated by $t_{11}^{(k)}$,
$t_{21}^{(k)}$, $k=1, \ldots, p_1$, $t_{22}^{(r)}$, $r=1, \ldots,
p_2$ and $t_{12}^{(m)}$, $m=p_2-p_1+1, \ldots, p_2$.

We will denote by $\bar{t}_{11}^{(k)}$, $\bar{t}_{21}^{(k)}$,
$\bar{t}_{22}^{(k)}$, $\bar{t}_{12}^{(k)}$ the images of the
generators of $W(\pi)$ in the graded algebra $\bar{W}(\pi)$.

Let $$T_{11}(u)=\sum_{i=0}^{p_1}t_{11}^{(i)}u^{p_1-i}, \qquad
T_{22}(u)=\sum_{i=0}^{p_2}t_{22}^{(i)}u^{p_2-i},$$

$$T_{21}(u)=\sum_{i=1}^{p_1}t_{21}^{(i)}u^{p_1-i}, \qquad
 T_{12}(u)=
\sum_{i=1}^{p_1}t_{12}^{(s+i)}u^{p_1-i}$$ and \ben
D_1(u)=T_{11}(u),\qquad
D_2(u)=T_{11}(u+1)T_{22}(u)-T_{21}(u+1)T_{12}(u). \een The
coefficients $d_1^{(1)}, \ldots, d_{p_1}^{(1)}$ of $D_1(u)$ and
$d_{1}^{(2)}, \ldots, d_{p_1+p_2}^{(2)}$ of $D_2(u)$ are
generators of the Gelfand-Tsetlin subalgebra $\Ga$. Denote by
$\bar{d}_i^{(j)}$ their images in the graded algebra.

Recall that a sequence ${x_1, \ldots, x_n}$ of elements of some
commutative ring $R$ is called  {\em regular} if, for all $i=1,
\ldots, n,$ the multiplication by $x_i$ is injective on
$$R/<x_1,\ldots, x_{i-1}>R$$ and $R/<x_1, \ldots, x_n>R\neq 0$.

Since $W(\pi)$ is a \emph{special filtered} algebra in the sense
of \cite{fo3}, by \cite[Theorem~1.1]{fo3} we only need to show
that $\bar{d}_1^{(1)}, \ldots, \bar{d}_{p_1}^{(1)}$,
$\bar{d}_{1}^{(2)}, \ldots, \bar{d}_{p_1+p_2}^{(2)}$ is a regular
sequence in $\bar{W}(\pi)$.

We will use the following standard result.
\ble\label{lemma-regular-zeros} A sequence of the form
$x_1,\ldots, x_r, y_1, \ldots y_t$, where $y_1, \ldots, y_t$ are
homogeneous elements of $A=\Bbbk[x_1, \ldots, x_q]$, $q>r$, is
regular in $A$ if and only if the sequence $\tilde{y}_1, \ldots,
\tilde{y}_t$ is regular in $\Bbbk[x_{r+1}, \ldots, x_q]$, where
$\tilde{y}_i(x_{r+1}, \ldots, x_q)=y_i(0,\ldots, 0, x_{r+1},
\ldots, x_q)$. \ele

Applying Lemma~\ref{lemma-regular-zeros} we reduce the problem to
the problem of regularity of the sequence of images
 of $\bar{d}_{1}^{(2)}, \ldots, \bar{d}_{p_1+p_2}^{(2)}$
 in $\bar{W}(\pi)/(D_1(u))$. Consider the
 first $p_2$ elements in this sequence.
Then the image of $\bar{d}_i^{(2)}$ coincides with
$\bar{t}_{22}^{(i)}$, $i=1, \ldots, p_2$. Hence, applying again
Lemma~\ref{lemma-regular-zeros} we reduce the problem to the
regularity of the sequence of images of $\bar{d}_{p_2+1}^{(2)},
\ldots, \bar{d}_{p_1+p_2}^{(2)}$ in
$\bar{W}(\pi)/(T_{11}(u),T_{22}(u))$. Denote these elements by
$z_{p_2+1}, \ldots, z_{p_1+p_2}$.

Consider the restricted Yangian $\Y_{p_2}(\gl_2)$ of level $p_2$
(see \cite{fmo}) generated by the coefficients
of the polynomials
$$
T_{11}'(u)=\sum_{i=0}^{p_2}t_{11}^{(i)}u^{p_2-i}, \qquad
T_{22}'(u)=\sum_{i=0}^{p_2}t_{22}^{(i)}u^{p_2-i},
$$
$$T_{21}'(u)=\sum_{i=1}^{p_2}t_{21}^{(i)}u^{p_2-i},\qquad
T_{12}'(u)=\sum_{i=1}^{p_2}t_{12}^{(i)}u^{p_2-i}.
$$
Let $D(u)'= T_{21}'(u+1)T_{12}'(u)$. Let $y_1, \ldots,
y_{p_1+p_2}$ be the graded images of the coefficients of $D(u)'$
in ${\Y}_{p_2}(\gl_2)/(T_{11}'(u), T_{22}'(u))$. The Yangian
$\Y_{p_2}(\gl_2)$ is free over its Gelfand-Tsetlin subalgebra
generated by $T_{11}(u)'$ and $T_{11}'(u+1)T_{22}'(u)'-D(u)'$ by
 \cite[Theorem~3.4]{fmo}. Since the sequence  $y_1, \ldots, y_{p_1+p_2}$
is obtained from a regular sequence in $\bar{\Y}_{p_2}(\gl_2)$ by
substituting zeros instead of some generators, then $y_1, \ldots,
y_{p_1+p_2}$ is regular by Lemma~\ref{lemma-regular-zeros}. Hence
its subsequence $y_{p_2+1}, \ldots, y_{p_1+p_2}$ is also regular.
Thus the variety $V(y_{p_2+1}, \ldots, y_{p_1+p_2})\subset
\Bbbk^{2p_2}$ is equidimensional. Now project this variety on
 the subspace $\Bbbk^{2p_1}$ by substituting zeros instead of $t_{12}^{i}, i=1, \ldots, s$ and $t_{21}^{i},
i=p_1+1, \ldots, p_2$. The resulting variety is again
equidimensional of pure dimension $p_1$. Moreover, this variety
coincides with the variety  $V(z_{p_2+1}, \ldots, z_{p_1+p_2})$
and therefore the sequence $z_{p_2+1}, \ldots, z_{p_1+p_2}$ is
regular.

Hence we proved

\bth\label{theorem-freeness-gl_2} $W(\pi)$ is free as a right
(left) module over the Gelfand-Tsetlin subalgebra.

\eth

Consider the following analog of the \emph{Kostant-Wallach map}
(\cite{KoW})
$$KW: \Specm \bar{W}(\pi)\simeq \Bbbk^{3p_1+p_2}\rightarrow \Specm \bar{\Ga}\simeq \Bbbk^{2p_1+p_2}.$$
In particular we showed

\bco\label{corol-equidimensional} The map $KW$ is surjective and
the variety $KW^{-1}(0)$ is equidimensional of pure dimension
$p_1$.

\eco

We also get a good estimate of the size of the fiber for any
$\bm\in \Sp$.

\begin{thm} Let $(p_1, p_2)$ are the rows of $\pi$. For any
$\bm\in \Specm \Ga$ the fiber of $\bm$ constists of at most $p_1!$
isomorphism classes of irreducible Gelfand-Tsetlin
$W(\pi)$-modules. Moreover, the dimension of the subspace of
$\bm$-nilpotents in any such module is bounded by $p_1!$.
\end{thm}

\begin{proof}
Since $W(\pi)$ is free over $\Ga$ and $\Ga$ is a polynomial ring
then all conditions of \cite[Theorem~5.3,(iii)]{fo-Ga2} are
satisfied, and hence we have that  the fiber of $\bm$ consists of
at most $p_1!(p_1+p_2)!$ isomorphism classes of irreducible
Gelfand-Tsetlin $W(\pi)$-modules. But this bound can be improved
following \cite[Corollary~6.1,(2)]{fo-Ga2}. Let $\bar{\Ga}$ be the
integral closure of $\Gamma$ in $L$. If $\ell\in \Specm \bar{\Ga}$
projects to $\bm\in \Specm \Gamma$ then we  write
$\ell=\ell_{\bm}$. Note that given $\bm\in \Specm \Ga$ the number
of different $\ell_{\bm}$ is finite. Moreover, for any $\bm\in
\Specm \Ga$ and some fixed $\ell_{\bm}$ there exists at most
$p_1!$ elements $s\in \mathbb{Z}^{p_1}$ such that $\ell_{\bm}$ and
$\ell_{\bm}+s$ differ by the action of $G=S_{p_1}\times
S_{p_1+p_2}$. It immediately implies the statement about the bound
for the fiber. The same number bounds the dimension of the
subspace of $\bm$-nilpotents by \cite[Corollary~6.1,(1)]{fo-Ga2}.
\end{proof}

\section*{Acknowledgment}

The authors acknowledge the support of the Australian Research
Council. The first author  is supported in part by the CNPq grant
(processo 301743/2007-0) and by the Fapesp grant (processo
2005/60337-2). The first author is grateful to the University of
Sydney for support and hospitality. The authors are grateful to
T.Levasseur and P.Etingof for helpful discussions on Noether
problem.

\end{document}